\RequirePackage{fix-cm}
\documentclass[smallextended,envcountsect,final=]{svjour3} 
\smartqed 
\usepackage{graphicx}
\usepackage{amssymb}
\usepackage{amsmath}
\usepackage{algorithm}
\usepackage{algpseudocode}	
\usepackage{algorithmicx}
\usepackage{mathrsfs}
\usepackage{bm}
\usepackage{xcolor}
\usepackage{pgfplots}
\usepackage{subfigure}
\usepackage{booktabs}
\usepackage{multirow}
\usepackage{url}
\usepackage{graphicx} 
\usepackage{setspace}

\journalname{Noname}

\newcommand{\change}[1]{\textcolor{black}{#1}}
\newcommand{\StartChange}{\color{black}}
\newcommand{\EndChange}{\color{black}}

\newcommand{\Proof}{\noindent {\it Proof }$\;\;$}
\newcommand{\tr}{\top}

\begin{document}

\title{Decentralized Optimization over Tree Graphs}

%\subtitle{Using  the  LaTex Template}

\author{Yuning Jiang \and Dimitris Kouzoupis \and Haoyu Yin \and Moritz Diehl \and Boris Houska$^*$\thanks{$^*$Corresponding author. %\email{borish@shanghaitech.edu.cn}
} }

\institute{
Y. Jiang, H. Yin and B. Houska \at
School of Information Science and Technology\\
ShanghaiTech University, Shanghai, China\\
\email{jiangyn, yinhy, borish@shanghaitech.edu.cn}\\[0.2cm]
D. Kouzoupis and M. Diehl \at
Department of Microsystems Engineering (IMTEK)\\
University of Freiburg, Germany\\
\email{dimitris.kouzoupis, moritz.diehl@imtek.uni-freiburg.de}
}

\date{Received: date / Accepted: date}
%The correct dates will be entered by the editor.

\maketitle

\begin{abstract}
This paper presents a decentralized algorithm for non-convex optimization over tree-structured networks. We assume that each node of this network can solve small-scale optimization problems and communicate approximate value functions with its neighbors based on a novel multi-sweep communication protocol. In contrast to existing parallelizable optimization algorithms for non-convex optimization the nodes of the network are neither
synchronized nor assign any central entity. None of the nodes needs to know the whole topology of the network, but all nodes know that the network is tree-structured. We discuss conditions under which locally quadratic convergence rates can be achieved. The method is illustrated by running the decentralized asynchronous multi-sweep protocol on a radial AC power network case study.
\end{abstract}
%\subclass{49J53 \and  49K99 \and more}

%All acknowledgements should be placed in the back of the paper after Conclusions..

%%%%%%%%%%%%%%%%%%%%%%%%%%%%%%%%%%%%%%%%%%%%%%%%%%%%%%
\section{Introduction}
Large-scale optimization problems over sparse networks arise in many applications such as resource allocation problems~\cite{Nedic2018a}, smart grid control problems~\cite{Braun2016},
traffic coordination problems~\cite{Hult2016,Jiang2017a},
power system operation problems~\cite{Engelmann2018,Molzahn2017}, and statistical learning problems~\cite{Boyd2011}. In terms of existing numerical algorithms,
which can solve such large-scale optimization problems, one typically distinguishes
between distributed and decentralized methods~\cite{Bertsekas1989},
as reviewed below.

\bigskip
\noindent
Distributed optimization algorithms typically parallelize most of their operations,
but communicate results to a central coordinator. This usually requires
one to synchronize the whole network, as the central operations can often only be
performed after receiving all results from the agents of the network. For example,
many variants of the Alternating Direction Method of Multipliers
(ADMM)~\cite{Boyd2011,Bertsekas1989,Shi2014} as well as the
Augmented Lagrangian Alternating Direction Inexact Newton (ALADIN) method~\cite{Houska2016} alternate between
solving small-scale optimization problems, which can be done in parallel by the agents of a network,
and solving large scale linear equation systems, which is typically done by a central coordinator. Here, the central solver runs
advanced sparse linear algebra routines, which can, optionally, parallelize part of their operations~\cite{Gondzio2009,Pakazad2017,Zavala2008}.
A complete overview of existing distributed optimization algorithms would go beyond
the scope of this paper, but an overview of distributed optimization methods, with a particular
focus on augmented Lagrangian approaches, can be found in~\cite{Hamdi2005,Hamdi2011,Bertsekas1979}.
Notice that modern distributed optimization algorithms are applicable to convex as well as
non-convex optimization problems. For example, convergence conditions for ADMM to local minimizers
of non-convex optimization problems can be found in~\cite{Hong2016}. Similarly, ALADIN has been
designed for solving non-convex optimization problems and conditions for global convergence
to local minimizers can be found in~\cite{Houska2016}.
In the context of the current paper, we also mention that there exist distributed variants of interior point methods
for convex optimization that have been tailored for tree sparse graphs, as analyzed in~\cite{Khoshfetrat2017}.

\bigskip
\noindent
In contrast to the above reviewed methods for distributed optimization, decentralized optimization methods
do not require a central coordination step. In the most general case, these algorithms are not synchronized
and the agents might not even know the whole network structure. However, one usually requires that
all agents can communicate with their neighbors~\cite{Bertsekas1989}. Classical decentralized optimization
methods are often based on dual decomposition~\cite{Bertsekas2014,Terelius2011}, but there also exist
decentralized consensus variants of ADMM~\cite{Shi2014,Makhdoumi2017}. However, in contrast to
distributed optimization methods, decentralized optimization algorithms often converge for convex problems
only and have a linear or even sublinear convergence behavior~\cite{Boyd2011,Nesterov2013,Nedic2018a}.
For instance, linear convergence conditions of a fully decentralized ADMM method for consensus optimization over networks
has been established in~\cite{Shi2014}.
A unified convergence analysis for decentralized ADMM in dependence on the network structure can also be
found in~\cite{Makhdoumi2017}.

\bigskip
\noindent
In summary, there has been a huge amount of research on distributed optimization algorithms, but there are, at the current status of research, no generic asynchronous decentralized optimization algorithms available that are $i)$ applicable to large-scale non-convex optimization problems and $ii)$ locally equivalent to Newton-type methods such that locally superlinear or quadratic convergence rates can be expected. Therefore, this paper asks the question whether such decentralized algorithms can be constructed at all, at least for special classes of networks. Here, our focus is on networks with tree graphs, which arise in some (but not all) applications such as traditional optimal control problems~\cite{Bellman1966} or receding horizon control problems~\cite{Rawlings2017}, where linear trees occur, scenario multi-stage MPC problems~\cite{Bernardini2011,Kouzoupis2018,Lucia2014}, which have less trivial tree structures, or radial power grid networks that possess non-trivial tree structure too~\cite{Kekatos2012,Peng2014}.

\bigskip
\noindent
Because there exist dynamic programming (DP) methods~\cite{Bertsekas2007} as well as approximate DP methods~\cite{Bertsekas2005,Keshavarz2014,Wang2015}, which can exploit the structure of optimization problems over tree-topologies, we briefly review these methods as well as their closely related min-sum algorithms~\cite{Kellerer2014,Kellerer2016} in Section~\ref{sec::TreeOpt}. The main contribution of this paper is then presented in Section~\ref{sec::multisweep}, where we develop a multi-sweep algorithm for solving tree-structured optimization problems. Moreover, Section~\ref{sec::multisweep2} presents a fully decentralized and asynchronous variant of this novel multi-sweep method, which can still achieve locally quadratic convergence rates. These theoretical developments are illustrated by a case study in Section~\ref{sec::caseStudies}, for a non-trivial radial power grid optimization problem.

\StartChange

\subsection*{Notation and Preliminaries}
Throughout this paper a couple of existing results from the field of parametric nonlinear optimization are used~\cite{Nocedal2006}. In order to briefly review these existing results, we consider a general (twice continuously differentiable) parametric equality constrained optimization problem of the form
\begin{eqnarray}
\label{eq::PARopt}
F(x) = \min_{z} f(x,z) \quad \text{s.t.} \quad c(x,z) = 0 \; \mid \; \lambda
\end{eqnarray}
Here, $f: \mathbb R^{n} \times \mathbb R^m \to \mathbb R$ and $c: \mathbb R^{n} \times \mathbb R^m \to \mathbb R^{n_c}$ denote the twice continuously differentiable objective and constraint functions, $x \in \mathbb R^n$ a parameter, and $z \in \mathbb R^m$ the primal optimization variable. Notice that we use the syntax ``$\mid \; \lambda$'' after an equality constraint to say that $\lambda$ denotes the multiplier of this constraint. Consequently, in the context of~\eqref{eq::PARopt}, we have $\lambda \in \mathbb R^{n_c}$.

\begin{definition}
\label{def::licq}
We say that a (local) minimizer $z^\star(x)$ of~\eqref{eq::PARopt} satisfies the linear independence constraint qualification (LICQ) at a given point $x \in \mathbb R^n$, if the constraint Jacobian matrix
\[
\frac{\partial c(x,z^\star(x))}{\partial z}
\]
has full-rank.
\end{definition}

\noindent
At this point, we recall that the LICQ condition is sufficient to ensure that minimizers of~\eqref{eq::PARopt} are KKT points. That is, if $z^\star(x)$ is a local minimizer of~\eqref{eq::PARopt} at $x$ at which LICQ holds, then there exists a unique multiplier $\lambda^\star(x) \in \mathbb R^{n_c}$ such that (see \cite[Thm.~12.1]{Nocedal2006})
\begin{eqnarray}
\label{eq::KKTA}
0 &=& \nabla_z L(x,z^\star(x),\lambda^\star(x)) \\[0.16cm]
\label{eq::KKTB}
0 &=& c(x,z^\star(x)) \; ,
\end{eqnarray}
where $L: \mathbb R^{n} \times \mathbb R^m \times \mathbb R^{n_c} \to \mathbb R$ denotes the Lagrangian function,
\[
L(x,z,\lambda) = f(x,z) + \lambda^\tr c(x,z) \; .
\]
In this context, another important regularity conditions is the so-called second order sufficient condition:

\begin{definition}
\label{def::sosc}
We say that a (local) primal-dual minimizer $(z^\star(x)$, $\lambda^\star(x))$ of~\eqref{eq::PARopt} satisfies the second order sufficient condition (SOSC) at $x$, if for any vector $v \in \mathbb R^n$ with $v \neq 0$ and
\[
\frac{\partial c(x,z^\star(x))}{\partial z} v = 0
\]
it follows that $v^\tr \nabla_z^2 L(x,z^\star(x),\lambda^\star(x)) v > 0$.
\end{definition}

\noindent
In some of the technical derivations of this paper, we assume that both LICQ and SOSC hold at a given point of an optimization problem. This motivates the following definition.

\begin{definition}
\label{def::reg}
We say that a local minimizer $z^\star(x)$ of~\eqref{eq::PARopt}, together with the multiplier $\lambda^\star(x)$,  is a regular minimizer at $x$, if it satisfies the LICQ and the SOSC conditions.
\end{definition}

\noindent
The above definition is useful, because it allows us to establish the following regularity result for the value function $F$, which follows by applying the implicit function theorem~\cite[Thm.~A2]{Nocedal2006} (or its generalized version~\cite{Robinson1980}) to the first order necessary KKT conditions~\eqref{eq::KKTA}-\eqref{eq::KKTB}; see also~\cite{Houska2013} for details.

\begin{lemma}
\label{lem::A1}
If $(z^\star(x),\lambda^\star(x))$ denotes a regular minimizer and multiplier of~\eqref{eq::PARopt} at~$x$, then $F$ is twice continuously differentiable at $x$ and we have
\[
\nabla_x F(x) = L_x(x) \qquad \text{and} \qquad \nabla_x^2 F(x) = L_{xx}(x) - L_{xu}(x) L_{uu}(x)^{-1} L_{ux}(x)
\]
with shorthands $L_x(x) = \frac{\partial L}{\partial x}(x,z^\star(x),\lambda^\star(x))$ as well as
\begin{eqnarray}
L_{xx}(x) &=& \frac{\partial^2 L}{\partial x^2}(x,z^\star(x),\lambda^\star(x)) \; , \notag \\[0.2cm]
L_{ux}(x) &=& L_{xu}(x)^\tr = \left(
\begin{array}{c}
\frac{\partial^2 L}{\partial z \partial x}(x,z^\star(x),\lambda^\star(x)) \\[0.16cm]
\frac{\partial c}{\partial x}(x,z^\star(x))
\end{array}
\right) \; , \notag \\[0.16cm]
\text{and} \qquad L_{uu}(x) &=& \left(
\begin{array}{cc}
\frac{\partial^2 L}{\partial z^2}(x,z^\star(x),\lambda^\star(x)) & \frac{\partial c}{\partial z}(x,z^\star(x))^\tr \\[0.16cm]
\frac{\partial c}{\partial z}(x,z^\star(x)) & 0
\end{array}
\right) \; ,
\end{eqnarray}
where $L_{uu}(x)$ is invertible at $x$. Moreover, there exists an open neighborhood $\mathcal N \subseteq \mathbb R^n$ of $x \in \mathrm{int}(\mathcal N)$ such that $F$ is twice differentiable on $\mathcal N$ and such that $(z^\star(x),\lambda^\star(x))$ is a regular solution of~\eqref{eq::PARopt} at $y$ for all $y \in \mathcal N$. 
\end{lemma}

Notice that the statement of this lemma can be strengthened further under the additional assumption that the second derivatives of $f$ and $c$ are Lipschitz continuous; see~\cite{Houska2013,Robinson1980} for details.

\begin{corollary}
\label{cor::A1}
Let $(z^\star(x),\lambda^\star(x))$ be a regular minimizer and multiplier of~\eqref{eq::PARopt} at $x$. If the second derivatives of $f$ and $c$ are locally Lipschitz continuous, then there exists an open neighborhood $\mathcal N$ of $x^\star$ such that the second derivatives of $F$ are locally Lipschitz continuous on $\mathcal N$.
\end{corollary}

\noindent
Throughout this paper we construct algorithms for solving nonlinear equality constrained problems of the form
\begin{eqnarray}
\label{eq::FC}
\min_x \; F(x) \quad \text{s.t.} \quad C(x) = 0  \; | \: \kappa \; , 
\end{eqnarray}
where $C$ typically denotes a twice Lipschitz-continuously differentiable (in the applications of this paper even affine) consensus constraint. Let
\[
\Phi(x,y) = F(y) + \mathbf{O}(\| x - y \|^3)
\]
denote a local model of $F$ at $x$ and consider an iteration of the form
\begin{eqnarray}
\label{eq::AlmostNewton}
x^{k+1} = \underset{y}{\text{argmin}} \; \Phi(x^k,y) \quad \text{s.t.} \quad C(y) = 0 \; | \: \kappa^{k+1}
\end{eqnarray}
started at an initial point $x_0 \in \mathbb R^n$. The following theorem is---at least in very similar versions---known and well-established~\cite[Thm.~18.4]{Nocedal2006}. However, for the sake of completeness, we provide a short proof in Appendix~\ref{app::proof}.

\begin{theorem}
\label{thm::AlmostNewton}
Let $x^\star$ be a regular minimizer of~\eqref{eq::FC}; that is, such that LICQ and SOSC are satisfied. If $F,C$ and $\Phi(x,\cdot)$ are twice differentiable functions with locally Lipschitz continuous second derivatives with uniform Lipschitz constant for all $x$ in a neighborhood of $x^\star$, then the iterates $x^k$ of~\eqref{eq::AlmostNewton} converge locally to $x^\star$ with quadratic convergence rate; that is,
\[
\| x^{k+1} - x^\star \| \leq \mathbf{O} \left( \Vert x^k - x^\star \Vert^2 \right)
\]
for all $k \in \mathbb N$ whenever the initialization $x^0$ is in a sufficiently small neighborhood of $x^\star$.
\end{theorem}

\EndChange

%%%%%%%%%%%%%%%%%%%%%%%%%%%%%%%%%%%%%%%%%%%%%%%%%%%%%%
\section{Optimization over tree graphs}
\label{sec::TreeOpt}
This section introduces optimization problems over tree graphs and briefly discusses the advantages and disadvantages of existing dynamic programming methods, which can be used to solve them numerically.

\subsection{Tree graphs}
Let $(\mathcal{N},\mathcal{E})$ denote a graph with node set $\mathcal{N}=\{1,...,N\}$
and edge set $\mathcal{E}\subseteq \mathcal{N}\times \mathcal{N}$. In the following, we assume that this graph is undirected; that is, $\mathcal{E}$ is symmetric: $(i,j) \in \mathcal E$
if and only if $(j,i) \in \mathcal E$. We call $(\mathcal{N},\mathcal{E})$
a tree graph if there exists for every $i \in \mathcal N \setminus \{ 1 \}$ exactly one way to walk from the first
node to the $i$-th node via the edges of the graph without passing any node more than once. Here, the first node is called
the root of the tree graph. However, it is important to keep in mind that the definition of a tree graph does
not depend on how the nodes are enumerated. In particular, any node in a tree can be a root as long as we
re-enumerate the nodes accordingly. An example for a simple tree structure is shown in the sketch below.
\begin{figure}[H]
	\centering
	\includegraphics[width=0.4\linewidth]{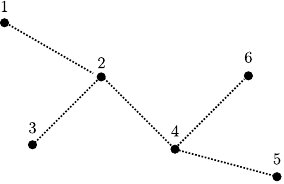}
	\caption{Example for a tree graph with $N=6$ nodes.}
	\label{fig::topology1}
\end{figure}
\noindent
In the following, we denote by $\mathcal N_i = \{ j \in \mathcal N \mid (i,j) \in \mathcal E \}$
the set of neighbors and by $\mathcal L = \{ i \in \mathcal N \mid |\mathcal N_i| \leq 1 \}$ the
set of nodes with at most one neighbor. For example, for the network in Figure~\ref{fig::topology1},
we have $\mathcal{L} = \{1,3,5,6\}$. In the following, we additionally use the notation
\[
\mathcal{L}^\bullet = \mathcal{L} \setminus  \{ 1 \}  \; ,
\]
which is called the set of leaves. Notice that the definition of $\mathcal L$ does not depend
on which node is assigned as root. This is in contrast to the set $\mathcal L^\bullet$ of leaves,
which may not contain the root node.
\begin{proposition}
	\label{prop::EP}
	The nodes of a tree graph $(\mathcal N, \mathcal E)$ can be enumerated in such a way that
	the graph $(\mathcal N, \mathcal E^+)$, with
	$$\mathcal E^+ = \{ (i,j) \in \mathcal E \mid i < j \} \; ,$$
	is still a tree graph.
\end{proposition}
The proof of the above proposition is constructive: we can start at any node and call it the root
by assigning the label $1$. Next, the root enumerates its children in increasing order, the
children enumerate their children, and so on, until all nodes have a number and
such that the set
\[
\mathcal C_i = \{ j \in \mathcal N_{i} \mid j \geq i \}
\]
corresponds to the set of children of the $i$-th node. Notice that this enumeration
procedure also ensures that the equation
\[
\mathcal{L}^\bullet = \{ i \in \mathcal N \mid |\mathcal C_i| = 0 \}
\]
holds; that is, leaves have no children. Last but not least, every node with number
$i \geq 2$ has a unique parent node $\pi_i \in \mathcal{N}_i\backslash\mathcal{C}_i$.

\subsection{Tree-structured optimization problems}

This paper concerns structured optimization problems of the form
\begin{equation}
\label{eq::opt}
V^* = \min_{x} \; \sum_{ i \in \mathcal{N}} F_{i}(x_i) \quad \text{s.t.} \quad \left\{
\begin{array}{l}
\forall (i,j) \in \mathcal E^+, \\[0.1cm]
S_{i,j}x_{i}=S_{j,i}x_{j}\; .
\end{array}
\right. 
\end{equation}
Here, the functions $F_i: \mathbb R^{n_i} \to \mathbb R$ denote objective functions,
$S_{i,j} \in \mathbb{R}^{n_{i,j}\times n_i}$ and $S_{j,i}\in \mathbb{R}^{n_{i,j}\times n_j}$ given connectivity matrices
and $(\mathcal N, \mathcal E)$ a tree graph, where $\mathcal E^+$ is defined as in Proposition~\ref{prop::EP}. Notice that in the most general case, there are no further assumptions on the functions $F_i$ needed as long as the minimizer of~\eqref{eq::opt} exists\change{, although for some of the algorithmic developments below, we will work with stronger local regularity assumptions on the functions $F_i$, such as twice Lipschitz-continuous differentiability. Notice that in the context of~\eqref{eq::opt} the function  $F_i$ could denote the minimum value of a parametric optimization problem~\cite{Kouzoupis2016b,Kouzoupis2019b} that is solved by
the $i$-th node in order to evaluate $F_i$ at the parameter~$x_i$,
\begin{eqnarray}
F_i(x_i) = \min_{z_i} f_i(x_i,z_i) \quad \text{s.t.} \quad c_i(x_i,z_i) \; ,
\end{eqnarray}
where the function $c_i$ can be used to enforce nonlinear equality constraints. As long as the regularity assumptions from Lemma~\ref{lem::A1} and Corollary~\ref{cor::A1} are satisfied for the functions $f_i$ and $c_i$, these results can be used to ensure that $F_i$ is locally twice continuously differentiable.}

Now, the goal of this paper is to develop decentralized optimization algorithms for solving~\eqref{eq::opt} allowing neighbor-to-neighbor communication only. Here, we are particularly interested in algorithms that specify decentralized communication protocols that are independent of the particular structure of the tree $(\mathcal N, \mathcal E)$.
This means that none of nodes should be required to know the complete structure of the tree.

\subsection{Dynamic programming}
Problem~\eqref{eq::opt} can be solved by implementing the \change{dynamic programming
recursion~\cite{Bellman1966,Bertsekas2013}}
\begin{equation}
\label{eq::dyn}
V_i(p)= \min_{x_i} \; F_i(x_i) + \sum_{j \in \mathcal C_i} V_{j}( S_{i,j} x_i ) \; \; \mathrm{s.t.} \; \;
S_{i,\pi_i}x = p
\end{equation}
for all $i\neq 1$. Notice that this optimization problem can be solved by the $i$-th node
as soon as it has received the value functions $V_j$ from all its children. Moreover,
the recursion starts at the leaves $i \in \mathcal L^\bullet$, because these nodes
do not have children. In the last backward recursion step, the root node solves the optimization
problem
\begin{equation}
\label{eq::dynRoot}
V^* = \min_{x_1} \; F_1(x_1) + \sum_{ j \in \mathcal C_1} V_{j}(S_{1,j} x_1)
\end{equation}
finding the optimal value $V^*$ of~\eqref{eq::opt}. Finally, the root initializes a so-called forward sweep by sending linear combinations, $S_{1,j} x_1^\star$, to its children such that they can
find their solutions $x_j^\star$, send linear combinations of these solutions to their
children, and so on, until all nodes know the optimal solution. Notice that a desirable advantage of
dynamic programming is that the globally optimal solution is found. Moreover, the number of communication steps of dynamic programming is relatively low: we need only one backward sweep from the leaves to the root and one forward sweep from the root to the leaves. However, one disadvantage of dynamic programming is that one needs to construct
the functions $V_{i}$\change{, which is only possible with high numerical precision if further regularity assumptions on the functions $F_i$ are introduced~\cite{Gruene2002,Luss1990}.} 
If the dimension $n_{i,j}$ of the coupling variables
%If we have $n_{i,j} = 1$ or if the dimension of the coupling variables 
is small, this is not a big problem, but, in general, dynamic programming is affected by the so-called curse of dimensionality.

\begin{remark}
Variants of dynamic programming for general tree structured networks have been developed in~\cite{Kellerer2014} under the name \textit{min-sum algorithms}. An approximate variant of these min-sum algorithms for piecewise quadratic optimization problems can be found in~\cite{Kellerer2016}, which is related to the developments in this paper, although we consider a much more general class of non-convex optimization problems.
\end{remark}

\subsection{Run-time considerations}

In order to briefly discuss the run-time properties of dynamic programming, we introduce the following definition.
\begin{definition}
	\label{def::d}
	We denote with $d(i,j)$ the minimum number of edges over which one has to walk in order to get from any node $i \in \mathcal N$ to a node $j \in \mathcal N$.
\end{definition}
\begin{proposition}
	The function $d: \mathcal N \times \mathcal N \to \mathbb N$, as introduced in Definition~\ref{def::d}, is a metric on $(\mathcal N,\mathcal E)$. Moreover, $d$ is invariant under re-enumeration of the nodes.
\end{proposition}
\Proof
It is easy to check that $d(i,j) = 0$ if and only if $i=j$, $d(i,j) = d(j,i)$, and $d(i,k) \leq d(i,j) + d(j,k)$ for all $i,j,k \in \mathcal N$
and $d$ is indeed a metric. The additional invariance statement follows trivially from the fact that the definition of $d$ depends only on how the nodes are connected by edges---not on how they are enumerated.
\qed

\bigskip
\noindent
Now, if the dynamic programming method is run in synchronous mode; that is, such that every node executes at most one dynamic programming step per sampling time $\delta > 0$, then the root node $1$ receives all value functions after time $\delta D$, where
\[
D = \max_{j} \, d(1,j)
\]
denotes the depth of the enumerated tree. Because the forward sweep takes equally long, the total run-time of a synchronized dynamic programming method with sampling time $\delta$ is given by
\begin{equation}
\label{eq::TD}
T = 2 \delta D \; ,
\end{equation}
where $D$ denotes the depth of the tree. This implies in particular that the run-time of dynamic programming (DP) is not invariant under re-enumeration of the nodes, because the depth of a tree depends on which node is assigned as the root node $1$.
\begin{example}
Let us come back to 
the network from Figure~\ref{fig::topology1}, whose depth is given by
$$\max_{j} \, d(1,j) = 3 .$$
Equation~\eqref{eq::TD} implies that the synchronized DP run-time of this network is $6 \delta$. However, if we would re-enumerate the
nodes of this network, such that the node with label $2$ becomes the new root, the synchronized DP run-time improves to $4 \delta$.
\end{example}

\section{Multi-sweep method}
\label{sec::multisweep}
In order to avoid the curse of dimensionality of standard dynamic programming, we replace~\eqref{eq::dyn} by an approximate dynamic programming recursion. 
For this aim, we introduce parametric auxiliary optimization problems of the form
\begin{equation}
\label{eq::dynApprox}
\begin{array}{rccl}
\Omega_{i,k}(p) &=&\underset{x_i}{\min} & F_i(x_i) + \sum\limits_{j \in \mathcal N_i \setminus \{k\}} W_{j,i}( S_{i,j} x_i ) \\[0.16cm]
& &\mathrm{s.t.}& S_{i,k}x_i = p
\end{array}
\end{equation}
for any $k \in \mathcal N_i$ and recursively constructed approximation functions $W_{i,k} \approx \Omega_{i,k}$.
At this point, one can, in principle, admit all kinds of approximation functions, but the main limitation
is that the functions $W_{i,k}$ should be representable in a suitable storage format,
such that they can be sent over the network links with reasonable effort.
For example, one could construct quadratic approximations that can be stored in the form of Hessian matrix, gradient vector and a scalar. Notice that if we set $k = \pi_i$ in~\eqref{eq::dynApprox} and propagate
the approximation in backward mode from the leaves to the root, then this construction is such that
$$W_{i,\pi_i} \approx \Omega_{i,\pi_i} \approx V_i \; ;$$
that is, $W_{i,\pi_i}$ can be interpreted as an approximation of the function $V_i$.
Notice that one can repeat the approximate dynamic programming recursion multiple times in order
to refine the accuracy of the approximation. In order to develop such a multi-sweep variant, we assume that the $i$-th node has two modes: a backward and a forward mode. Moreover, the $i$-th node (with $i \neq 1$) is initially set
to backward mode and it is initialized with a model $W_{\pi_i,i}$ of its parent. If no such model is
available, we may, for example, set $W_{\pi_i,i} = 0$. Now, if the $i$-th node is in backward mode,
it waits until it receives models $W_{j,i}$ from all children $j \in \mathcal C_i$, and then solves
the auxiliary optimization problem
\begin{equation}
\label{eq::NLP}
\min_{y_i} \;\; F_i(y_i) + \sum_{j \in \mathcal N_i} W_{j,i}( S_{i,j} y_i ) \; .
\end{equation}
Let $y_i^\star$ denote a minimizer of this problem. Next, the approximation $W_{i,\pi_i} \approx \Omega_{i,\pi_i}$
in~\eqref{eq::dynApprox} is constructed in such a way that we have
\begin{equation}
\label{eq::ApproxOrder}
W_{i,\pi_i}(p) = \Omega_{i,\pi_i}( p ) + \mathbf{O}( \left\|  S_{i,\pi_i} y_i^\star - p \right\|^{q+1} ) \; ,
\end{equation}
where $q$ denotes the order of the approximation, and sent to the parent node $\pi_i$. After this, the $i$-th node
is set to forward mode. In this mode, it waits until a new model $W_{\pi_i,i}$ is sent from the parent node, then~\eqref{eq::NLP}
is solved once more in order to update $y^\star$. The corresponding updated models $W_{i,j} \approx \Omega_{i,j}$, with
\begin{equation}
\label{eq::ApproxOrderBack}
W_{i,j}(p) = \Omega_{i,j}( p ) + \mathbf{O}( \left\|  S_{i, j} y_i^\star - p \right\|^{q+1} ) \; ,
\end{equation}
are sent to all children $j \in \mathcal C_i$. Notice that the protocol for the root node, $i =1$, is completely
analogous with the only difference being that this node immediately switches to forward mode as soon as the models
from all children are received. The complete multi-sweep procedure is summarized in Algorithm~\ref{alg::MDP}.
\begin{algorithm}[h]
	\setstretch{1.2}
	\caption{Multi-sweep method of order $q$}
	\label{alg::MDP}
	\begin{algorithmic}[1]
		\State \textbf{Initialization:}
		\State Set all nodes to backward mode.
		\State If $i \neq 1$, choose an initial model function $W_{\pi_i,i}$.
		\State \textbf{Repeat the following protocol on every node $i \in \mathcal N$:}
		\State \textbf{IF} the node is in backward mode:
		\State Wait for model updates $W_{j,i}$ from all children $j \in \mathcal C_i$.
		\State If $i=1$, switch to forward mode and \textbf{BREAK}
		\State Solve~\eqref{eq::NLP} and construct $W_{i,\pi_i} \approx \Omega_{i,\pi_i}$ such that~\eqref{eq::ApproxOrder} holds.
		\State Send $W_{i,\pi_i}$ to the parent and switch to forward mode. 
		\State \textbf{ELSE} (the node is in forward mode):
		\State If $i \neq 1$, wait for the model update $W_{\pi_i,i}$ from the parent.
		\State Solve~\eqref{eq::NLP} and construct $W_{i,j} \approx \Omega_{i,j}$ such that~\eqref{eq::ApproxOrderBack} holds.
		\State Send $W_{i,j}$ to all nodes $j \in \mathcal C_i$, switch to backward mode.
		\State \textbf{ENDIF}
	\end{algorithmic}
\end{algorithm}

\subsection{Termination conditions}

Notice that once the protocols on all nodes are ``switched on'', Algorithm~\ref{alg::MDP} keeps on updating its model functions forwever. If one wishes to introduce a termination condition, this can be done by modifying Line~9 of Algorithm~\ref{alg::MDP} as follows: if the $j$-th node does not only send the model function $W_{j,\pi_j}$ but also the projected solution $S_{j,\pi_j} y_j^\star$ to its parent after solving~\eqref{eq::NLP}, then its parent can evaluate the residual
\[
r_i=\max_{j\in\mathcal{C}_i}\left\|S_{i,j}y_i^\star- S_{j,\pi_j}y_j^\star\right\|_{\infty} \; .
\]
Thus, if all these residual values are forwarded, the root node can evaluate the infinity norm, $r$, of the primal equality constraint violation as
\begin{equation}
\label{eq::stop}
r = \max_{i \in \mathcal N \setminus \mathcal L} \; r_i \; .
\end{equation}
If this constraint violation is small, $r \leq \epsilon$, for a small numerical tolerance $\epsilon > 0$, the root can send out a termination message that can be forwarded by the children until all nodes terminate. Notice that this termination condition merely ensures that the primal consensus constraint violation is small, but, in general, this is not sufficient to ensure that the solutions $y_i^\star$ of the nodes are close to a minimizer $x_i^\star$ of~\eqref{eq::opt} upon termination. Nevertheless, if one assumes that further regularity assumptions hold, for example if all minimizers of~\eqref{eq::opt} are regular KKT points, one can show that the termination condition $r \leq \epsilon$ ensures
\[
\Vert x^\star - y^\star \Vert \leq \mathbf{O}(\epsilon) ,
\]
see~\cite{Kouzoupis2019b,Nocedal2006} for details.

\subsection{Construction of model functions}
As mentioned in the previous section, Algorithm~\ref{alg::MDP} can, in the most general case, be applied
without further assumptions on the functions $F_i$ as long as one ensures that all minimizers are well-defined.
However, if one is interested in constructing practical
algorithms with $q \geq 1$, one might be interested in matching the first $q$ derivatives of the functions
$W_{i,k}$ and $\Omega_{i,k}$, which is only possible if the $F_i$s are sufficiently often differentiable and if the
minimizers of~\eqref{eq::dynApprox} are regular KKT points for all possible evaluation points $p$. A practical example is summarized below.

\begin{example}
\label{ex::cubicModel}
Let us assume that $F$ is three times continuously differentiable and has bounded third derivatives and that~\eqref{eq::dynApprox} has a regular parametric minimizer such that $\Omega_{i,k}$ is twice differentiable with bounded third order weak derivatives \change{(see Corollary~\ref{cor::A1})}. In this case, the models
$$W_{i,k}(p) = \frac{1}{2} p^\tr H_{i,k} p + g_{i,k}^\tr p + \sigma_{i,k} \Vert p-S_{i,k}y_i^\star \Vert^3  + \mathrm{const.}$$
can be constructed by setting
\[
H_{i,k} = \nabla^2\, \Omega_{i,k}( S_{i,k} y_i^\star)
\]
and
\[
g_{i,k} = \nabla\, \Omega_{i,k}( S_{i,k} y_i^\star) - H_{i,k} S_{i,k} y_i^\star \; .
\]
\change{Notice that these first and second order derivatives can be computed easily by using Lemma~\ref{lem::A1}.}
Here, $\sigma_{i,k} \geq 0$ is a cubic regularization constant, which is chosen such that $\sigma_{i,k} \geq \frac{1}{6}\Vert \nabla^3 \Omega_{i,k}(p) \Vert$ is a bound on the third order (weak) derivatives of $\Omega_{i,k}$ on a suitably defined validity domain of the model. Notice that a Taylor expansion of the function $\Omega_{i,k}$ at the point $S_{i,k} y_i^\star$ yields the relation
\begin{eqnarray}
\Omega_{i,k}(p) &\leq& \underbrace{\Omega_{i,k}(S_{i,k}y_i^\star) + y_i^\star S_{i,k}^\tr H_{i,k} S_{i,k} y_i^\star}_{\mathrm{const.}} \notag \\[0.16cm]
& & + g_{i,k}^\tr p + \frac{1}{2} p^\tr H_{i,k} p + \sigma_{i,k} \Vert p-S_{i,k}y_i^\star \Vert^3 \; .
\end{eqnarray}
Thus, our particular construction of $W_{i,k}$ is such that $W_{i,k}(p) \geq \Omega_{i,k}(p)$ as long as the constant offset of $W_{i,k}$ is chosen appropriately. Moreover, the associated multi-sweep method has order~$2$, since the above Taylor expansion based construction implies that
\[
W_{i,k}(p) = \Omega_{i,k}(p) +\mathbf{O}(\left\| p - S_{i,k} y_i^\star \right\|^3) \; .
\]
\end{example}
Notice that, in the above example, one only needs to send symmetric matrices $H_{i,k}$, vectors $g_{i,k}$ and $y_i^\star$, and the regularization constant $\sigma_{i,k}$ over the link $(i,k)$, as constant offsets of objective functions do not affect the optimal solutions of~\eqref{eq::NLP}. \change{The complexity for storing these three variables, including the matrix $H_{i,k}$, is given by $\mathbf{O}(n_{i,k}^2)$ recalling that $n_{i,k}$ denotes the number of coupling variables between the $i$-th and the $k$-th node of the network.}

\subsection{Conservation laws and convergence}
The convergence properties of Algorithm~\ref{alg::MDP} depend on the particular construction of the approximations $W_{i,k} \approx \Omega_{i,k}$
and their relation to the exact value functions $V_{i,k}$, which we define recursively as
\begin{equation}
\label{eq::DynAux}
\begin{array}{rccl}
V_{i,k}(p) &=&\underset{x_i}{\min} & F_i(x_i) + \sum\limits_{j \in \mathcal N_i \setminus \{k\}} V_{j,i}( S_{i,j} x_i ) \\[0.16cm]
& &\mathrm{s.t.}& S_{i,k}x = p
\end{array}
\end{equation}
for all $k \in \mathcal N_i$ and all $i \in \mathcal N$.
In this context, the following general conservation laws are useful for the composition of convergence conditions.

\begin{lemma}
\label{lem::bounds}
Algorithm~\ref{alg::MDP} has the following properties.
\begin{enumerate}
\item If all approximations $W_{i,k} \approx \Omega_{i,k}$ are lower bounds, such that $W_{i,k}(p) \leq \Omega_{i,k}(p)$ for all $p$, then
$$W_{i,k}(p) \leq V_{i,k}(p)$$
holds globally for all $p$, for all $k \in \mathcal N_i$ all $i \in \mathcal N$, and during all iterations.
\item The above statement also holds after replacing all $\leq$ signs with $\geq$ signs; that is, the conservation of upper bounds holds, too.
	\end{enumerate}
\end{lemma}

\Proof
For the case $|\mathcal N| = 1$ the statement of the lemma is trivial and we, thus, assume $|\mathcal N| > 1$.
The proof of the first statement of this lemma follows by a tree-structured induction that starts with $i \in \mathcal L^\bullet$ and propagates
through the tree. Here, our induction start uses that the inequality\footnote{The assumption $|\mathcal N| > 1$ ensures that $\pi_i$ exists and is well-defined for all $i \in \mathcal L^\bullet$.}
\[
W_{i,\pi_i}(p) \leq \Omega_{i,\pi_i}(p) = V_{i,\pi_i}(p)
\]
holds for all $i \in \mathcal L$ and their associated parent nodes $\pi_i \in \mathcal N_i$.
Next, our induction assumption is that
$$W_{k,i}(p) \leq V_{k,i}(p)$$
holds at a given node $i \in \mathcal N$ for all children $k \in \mathcal C_i$. Now, the definition
of $\Omega_{i,k}$ and $V_{i,k}$ in~\eqref{eq::dynApprox} and~\eqref{eq::DynAux} implies that
\begin{equation}\label{eq::UpperLowerBound}
\Omega_{i,\pi_i}(p) \leq V_{i,\pi_i}(p) \quad \Rightarrow \quad W_{i,\pi_i}(p) \leq V_{i,\pi_i}(p) \; ,
\end{equation}
where the latter statement can be interpreted as an intermediate induction conclusion yielding
that the inequality
$$W_{i,\pi_i}(p) \leq V_{i,\pi_i}(p)$$
holds for all $i \in \mathcal N \setminus \{ 1 \}$.
Similarly, the same induction argument can be repeated in forward mode, which yields that
\[
W_{i,k}(p) \leq V_{i,k}(p)
\]
also holds for all children $k \in \mathcal C_i$ and all $i \in \mathcal N \setminus \mathcal L^\bullet$.
The proof of the second statement of the lemma is completely analogous, as we can replace all
$\leq$ signs with $\geq$ signs without altering the logic of the proof.
\qed

\bigskip
\noindent
The above lemma can be used as a basis for convergence proofs of Algorithm~\ref{alg::MDP}. For example, a locally quadratic convergence statement for second order variants of Algorithm~\ref{alg::MDP} can be summarized as follows.

\begin{theorem}
\label{thm::local}
Let us assume that the approximation functions in Algorithm~\ref{alg::MDP} satisfy $W_{i,k}(p) \geq \Omega_{i,k}(p)$. \change{If the functions $F_i$ as well as the functions $W_{i,k}$ are all twice differentiable with Lipschitz continuous second derivatives,} if Algorithm~\ref{alg::MDP} has order $q = 2$ and if all leaves are initialized in a local neighborhood of a regular minimizer $x^\star$ of~\eqref{eq::opt}; that is, such that $y_i^\star$ is in a local neighborhood of $x_i^\star$ for all $i \in \mathcal L^\bullet$, then the iterates of Algorithm~\ref{alg::MDP} converge with locally quadratic convergence rate.
\end{theorem}

\Proof
\change{Because $x^\star$ is assumed to be regular minimizer of~\eqref{eq::opt} and because we assume that the functions $F_i$ are twice differentiable with Lipschitz continuous second derivatives, a recursive application of Corollary~\ref{cor::A1} to the dynamic programming recursion yields that the functions $V_{i,k}$ are locally twice differentiable with Lipschitz continuous second derivatives. Next, b}ecause the second statement in Lemma~\ref{lem::bounds} ensures that 
\begin{equation}\label{eq::Ineq1}
W_{i,k}(p) \geq \Omega_{i,k}(p)\quad\Rightarrow \quad W_{i,k}(p) \geq V_{i,k}(p)\;,	
\end{equation}
the iterates $y_i^\star$ are stable and remain in a neighborhood of $x_i^\star$. Moreover, because we assume that Algorithm~1 has order $q=2$, we have
\[
W_{i,k}(p) \leq \Omega_{i,k}(p) + \mathbf{O}(\| p - S_{i,k} y_i^\star\|^3)
\]
and Lemma~\ref{lem::bounds} can be used to propagate lower bounds, too, finding
\begin{equation}
\label{eq::Ineq2}
W_{i,k}( p ) \leq V_{i,k}( p ) + \mathbf{O}( \| p - S_{i,k} y_i^\star \|^3 ) \;.
\end{equation}
Thus, by using inequalities~\eqref{eq::Ineq1} and~\eqref{eq::Ineq2}, it follows that
\[
W_{i,k}(p) = V_{i,k}(p) + \mathbf{O}( \Vert p - S_{i,k} y_i^\star \Vert^3 )
\]
\change{for all $k \in \mathcal N_i$ all $i \in \mathcal N$.} Consequently, Algorithm~\ref{alg::MDP} is locally equivalent to the exact dynamic programming method~\cite{Bertsekas1979} up to terms of order $3$. \change{Moreover, the functions $V_{i,k}$ and $W_{i,k}$ are all locally twice differentiable with Lipschitz continuous second derivatives. Consequently, using once more that $x^\star$ is a regular minimizer, we can apply Theorem~\ref{thm::AlmostNewton} to show that the iterates converge with locally quadratic convergence rate.}
\qed
\begin{remark}
	Notice that the conditions of Theorem~\ref{thm::local} are satisfied for the approximation functions that have been constructed in Example~\ref{ex::cubicModel}. Moreover under the additional assumptions that the cubic regularization constants $\sigma_{i,k}$ are sufficiently large, the third derivatives of $F_i$ are bounded, and all KKT points of~\eqref{eq::NLP} are regular, one can show that Algorithm~\ref{alg::MDP} converges globally to stationary points of~\eqref{eq::opt}---this convergence result is obtained in analogy to Nesterov's cubic regularization method~\cite{Nesterov2006}.
\end{remark}

\section{Simultaneous multi-sweep method}
\label{sec::multisweep2}
This section develops a variant of Algorithm~\ref{alg::MDP} that is invariant with respect to permutations of the enumeration of the nodes of the tree graph. Here, the main idea is to start an approximate dynamic programming recursion at all nodes $i \in \mathcal L$ simultaneously without assigning a root recalling that the definition of $\mathcal L$ does not depend on the enumeration of the nodes. This leads to a simultaneous multi-sweep method as summarized in Algorithm~\ref{alg::SMM}.
\begin{algorithm}[H]
\setstretch{1.2}
\caption{Simultaneous multi-sweep method}
\label{alg::SMM}
\begin{algorithmic}[1]
\State \textbf{Initialization of node $i$:} Set $\mathcal R = \varnothing$ and $\ell = 0$.
\State \textbf{Repeat the following protocol on every  node $i \in \mathcal N$:}
\State If $\ell = 0$, denote with $\mathcal J \subseteq \mathcal N_i$ the indices of all neighbors that have sent updates $W_{j,i}$ after the last reset of the collection $\mathcal R$ and update
$$\mathcal R \leftarrow \mathcal R \cup \mathcal J \; .$$
If $|\mathcal R| = |\mathcal N_i|$ set $\ell = i$.
\State If $\ell > 0$ and $\ell \neq i$, wait until a model update $W_{\ell,i}$ from the $\ell$-th node is arriving.
\State If $\ell > 0$, solve~\eqref{eq::NLP} and construct $W_{i,j} \approx \Omega_{i,j}$ for all indices $j \in \mathcal N_i$
such that
\[
W_{i,j}(p) = \Omega_{i,j}(p) + \mathbf{O}( \Vert S_{i,j}y_i^\star - p \Vert^{q+1}) \; .
\]
If $\ell \neq i$, send $W_{i,j}$ to all neighbors $j \in \mathcal N_i \setminus \{ \ell \}$. Otherwise,
send $W_{i,j}$ to all neighbors $j \in \mathcal N$.
Set $\ell = 0$, and reset $\mathcal R = \varnothing$.
\State If $|\mathcal R| = |\mathcal N_i|-1$ with $\{ k \} = \mathcal N_i \setminus \mathcal R$,
solve~\eqref{eq::NLP} and construct $W_{i,k} \approx \Omega_{i,k}$ such that
\[
W_{i,k}(p) = \Omega_{i,k}(p) + \mathbf{O}( \Vert S_{i,k}y_i^\star - p \Vert^{q+1}) \; .
\]
If $k \neq i$, apply a short randomly chosen time delay: if the the line $(i,k)$ is not blocked, block the line $(k,i)$ and send $W_{i,k}$ to node $k$, set $\ell = k$ and unblock all lines again. Otherwise, if the line $(i,k)$ is blocked, skip.
\end{algorithmic}
\end{algorithm}
\noindent
Notice that the communication protocol for the nodes $i \in \mathcal N$ does not specify a-priori which neighbors act as children and which as parent. Here, every node keeps a local integer variable $\ell$, which is set to $0$ whenever the node is in backward mode. While being in backward mode, Line~3 allows this node to collect updated models from all of its neighbors. This step assumes that Node $i$ has a buffer such that local copies of all arriving models can be stored temporarily. The required maximum storage capacity of this buffer can be determined a-priori and for each node separately, as long as every node knows its number of neighbors. Notice that Step~3 assigns the value $\ell = i$, if it receives updates from all neighbors. In this case, the $i$-th node decides spontaneously to act as root. Moreover, Step~6 ensures that the node forwards an updated model to a spontaneously assigned parent as soon as it receives model updates from at least $|\mathcal N_i|-1$ neighbors. Last but not least, Step~5 ensures that the node performs updates and sends out new models to all children whenever it is switched to forward mode.  

\bigskip
\noindent
Because the protocol in Algorithm~\ref{alg::SMM} does not specify which node is assigned as root, it is possible that different nodes spontaneously act as a root while the algorithm is running. Thus, there arises the question how many nodes can act as a root at the same time and, more generally, what the precise differences between  Algorithm~\ref{alg::SMM} and Algorithm~\ref{alg::MDP} are. In order to give answers to these questions, we first analyze Algorithm~\ref{alg::SMM} under the additional assumption that all nodes act synchronously, based on a global clock and equal sampling time. However, a completely asynchronous implementation of Algorithm~\ref{alg::SMM} is possible, too, as discussed further below (see Section~\ref{sec::async}).
\begin{figure}[htbp!]
	\centering
	\includegraphics[width=\linewidth]{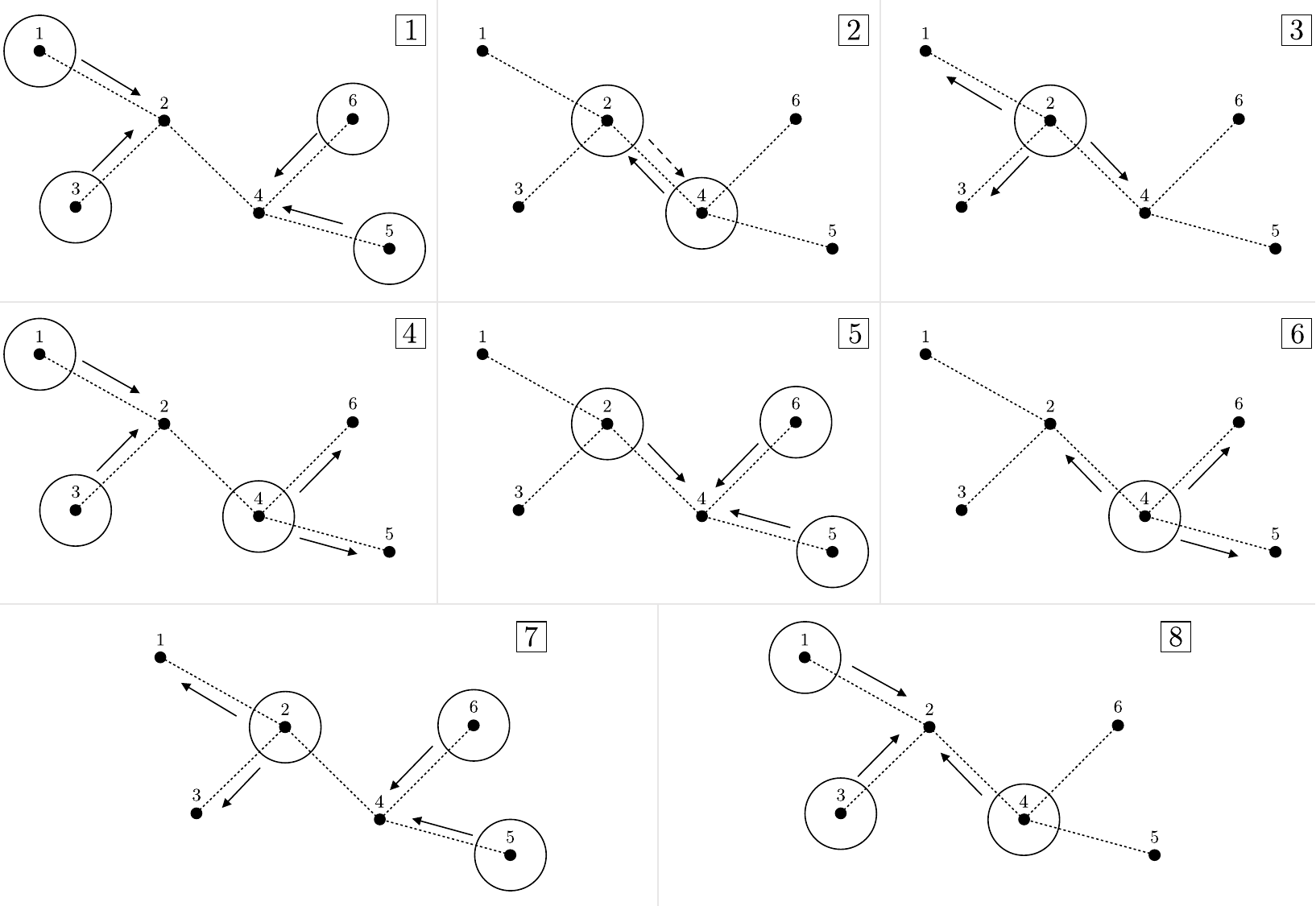}
	\caption{Visualization of the first $8$ iterations of Algorithm~\ref{alg::SMM} for a tutorial network with $6$ nodes.}
	\label{fig::example}
\end{figure}
\begin{example}
	Figure~\ref{fig::example} visualizes the first $8$ iterations that are obtained when executing Algorithm~\ref{alg::SMM} on all $6$ nodes of our tutorial network. In the first iteration, because we initialize with $\mathcal R = \varnothing$, the condition $|\mathcal R| = |\mathcal N_i| - 1$ in Step~6 of Algorithm~\ref{alg::SMM} is satisfied for all nodes $i \in \mathcal L = \{ 1,3,5,6\}$. Thus, these nodes solve optimization problems in parallel (indicated by circles) and send approximate value functions to the central nodes (indicated by arrows). In the second iteration, the nodes with labels $2$ and $4$ are solving their optimization problems. Here, we assume that Node~$4$ happens to be quicker and sends its result to Node~$2$ (indicated by a solid arrow), but blocks the reverse line (indicated by a dashed arrow). Thus, in the third iteration, Node $2$ acts spontaneously as root (see Step~5 of Algorithm~\ref{alg::SMM}) and sends its updated approximate value functions to all neighbors, $\mathcal N_2 = \{ 1,3,4\}$. Notice that the iterations continue with some nodes acting in forward mode while others are in backward mode. There is no preassigned root. For example, in the sixth iteration Node~$4$ happens to act has root, while in the third iteration the second node had taken this role.
	% If the nodes are not synchronized, in principle, any node can eventually become a root at some point in time.
\end{example}

\subsection{Synchronous multi-sweeps}
\label{sec::sync}
If the protocols of Algorithm~\ref{alg::SMM} are executed synchronously, the global behavior of this algorithm can be classified in dependence on the cardinality of the set of central nodes, which is defined as follows.
\begin{definition}
	Let $\Gamma \subseteq \mathbb N_+$ denote the set of central nodes,
	\[
	\Gamma = \underset{i \in \mathcal N}{\mathrm{\mathop{argmin}}} \, \max_{j \in \mathcal N} \, d(i,j) \; .
	\]
	We call $|\Gamma|$ the parity of the (undirected) graph $(\mathcal N, \mathcal E)$.
\end{definition}
Notice that the set $\Gamma$ can be interpreted as the set of nodes, which, if assigned as root, leads to a dynamic programming implementation with minimal run-time, since $\Gamma$ is the set of all nodes $i$, whose depth is minimal. For example, for the tree structured network in Figure~\ref{fig::topology1}, we have $\Gamma = \{2,4\}$ with parity $|\Gamma| = 2$.

\begin{lemma}
	\label{lem::parity}
	The parity of a tree structured graph is either even, $|\Gamma| = 2$, or odd, $|\Gamma| = 1$. Moreover, if it is even, then the two central nodes are neighbors.
\end{lemma}

\Proof
Let $L = \max\limits_{i,j} d(i,j)$ denote the length of a longest path in $(\mathcal N, \mathcal E)$ and
$$\mathbb L = \{ (i,j) \mid d(i,j) = L \}$$
the set of pairs $(i,j)$ with maximum distance. We distinguish two cases.

\smallskip
\noindent
\textbf{Case 1:} $L$ is even. In this case, there exists for every $(i,j) \in \mathbb L$ an odd number of nodes on the shortest path from $i$ and $j$, which implies that there is a (unique) central node $k$ on this path with $$d(i,k) = d(j,k) = \frac{L}{2} \; .$$
Let us assume that there is another pair $(i',j') \in \mathbb L$, whose central node $k'$ is not equal to $k$. Then, we must have
\begin{align}
d(i,k') + d(k',j) > d(i,j) = L \;,
\end{align}
since $k'$ is not on the shortest path from $i$ to $j$. But this means that either $d(i,k') > L/2$ or $d(j,k') > L/2$, which is a contradiction, as $d(i,k') \leq L/2$ and $d(j,k') \leq L/2$ must hold due to the construction of the central node $k'$. Thus, in summary, all pairs $(i,j) \in \mathbb L$ share the same central node $k$ and we must have $\Gamma = \{ k \}$ by construction.

\smallskip
\noindent
\textbf{Case 2:} $L$ is odd. In this case, there exist for every $(i,j) \in \mathbb L$ two central nodes $k_1,k_2$ with $d(k_1,k_2) = 1$ and such that
$$d(k_1,i) = (L-1)/2 \quad \text{and} \quad d(k_2,j) = (L+1)/2 \; .$$
Now, one can use a similar argument as in Case~1 to show that all $(i,j) \in \mathbb L$ share the same central nodes, finding that $\Gamma = \{ k_1, k_2 \}$.\\
\noindent
Both cases together yield the statement of this lemma.
\qed

\bigskip
\noindent
An immediate consequence of the above lemma is that the synchronized version of Algorithm~\ref{alg::SMM} is actually equivalent to Algorithm~\ref{alg::MDP} with one of the central nodes acting as root. For the case that the graph is odd, $|\Gamma| = 1$, the central root node is unique. Otherwise, for $|\Gamma| = 2$, one of the central nodes acts as root during every complete backward-forward sweep. However, in general, it cannot be predicted a-priori which of these central nodes act as root, as we have introduced the short random time delay in Step 6 of the protocol in Algorithm~\ref{alg::SMM}.
\begin{corollary}
	Algorithm~\ref{alg::SMM} is equivalent to a variant of Algorithm~\ref{alg::MDP}, where one of the central nodes acts as root during one complete backward-forward sweep. In particular, Algorithm~\ref{alg::SMM} converges under the same assumptions as~Algorithm~\ref{alg::MDP}.
\end{corollary}

\Proof
Notice that the proof of Lemma~\ref{lem::bounds} uses an induction argument under the assumption that the root node is fixed. In Algorithm~\ref{alg::SMM} the root is assigned dynamically, but there is one unique root during each backward-forward sweep, which means that the induction argument from the Lemma~\ref{lem::bounds} remains valid during every such complete sweep. But this means that the implications
\begin{eqnarray}
W_{i,k}(p) &\geq& \Omega_{i,k}(p)\quad\Rightarrow \quad W_{i,k}(p) \geq V_{i,k}(p) \notag \\[0.16cm]
\text{and} \qquad W_{i,k}(p) &\leq& \Omega_{i,k}(p)\quad\Rightarrow \quad W_{i,k}(p) \leq V_{i,k}(p)
\end{eqnarray}
also hold for Algorithm~\ref{alg::SMM} observing that these relations are independent of which node is assigned as root. By using once more that Algorithm~\ref{alg::SMM} has a unique root during each sweep, the local convergence rate estimate argument from the proof of Theorem~\ref{thm::local} remains valid, too. Thus, the statement of this corollary is a direct consequence of Theorem~\ref{thm::local}.
\qed

\subsection{Asynchronous multi-sweeps}
\label{sec::async}
An important observation of the previous section is that the root node is not assigned a-priori, but online while the algorithm is running. This allows us to run Algorithm~\ref{alg::SMM} in asynchronous mode. In this case, every node executes the protocol from Algorithm~\ref{alg::SMM} repeatedly, without coordinating the sampling time with other nodes. The analysis of this asynchronous variant is basically analogous to the synchronous case, but any node can act as the root node---not only the central nodes. The convergence analysis is, however, unaffected; that is, running Algorithm~\ref{alg::SMM} in asynchronous mode has no disadvantages in terms of its convergence properties.

\section{Application to a radial AC power network}
\label{sec::caseStudies}
\begin{figure}[H]
	\centering	
	\includegraphics[width = 0.95\linewidth]{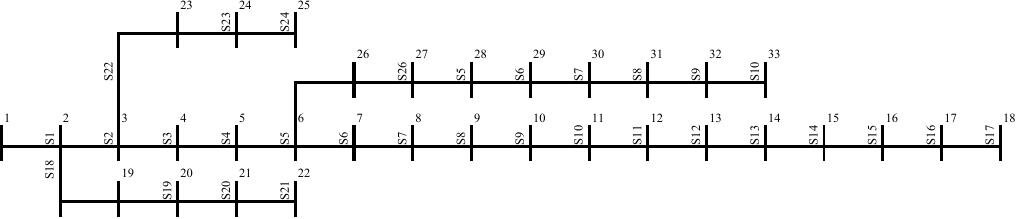}
	\caption{Topology of IEEE $33$-bus benchmark.}
	\label{fig::33bus}
\end{figure}

This section applies Algorithm~\ref{alg::SMM} to a state estimation problem for a radial AC power network. We use \texttt{MATPOWER v7.0}~\cite{Zimmerman2011} to generate an IEEE $33$-bus benchmark radial AC power network as shown in Figure~\ref{fig::33bus}.
\noindent
For this network we have  $\mathcal N = \{ 1, \ldots, 33 \}$ and $\mathcal E$ is defined precisely as visualized in the figure. The following section briefly reviews the model equations for power networks~\cite{Engelmann2018} and the related least-squares state estimation problem~\cite{Du2019}, which is here used as a benchmark problem for decentralized optimization.

\subsection{Model equations}
Let $v_i$ denote the voltage magnitude and $\theta_i$ the phase shift at the $i$-th node of the power network. The active and reactive power at this node is given by the power flow equations~\cite{Engelmann2018},
\begin{align}
p_i(v,\theta) &= v_{i} \sum_{j\in\mathcal{N}_i}v_{j}(G_{ij}\cos(\theta_i-\theta_j)+B_{ij}\sin(\theta_i-\theta_j))\notag \\
q_i(v,\theta) &= v_{i} \sum_{j\in\mathcal{N}_i}v_{j}(G_{ij}\sin(\theta_i-\theta_j)+B_{ij}\cos(\theta_i-\theta_j)) \; .\notag
\end{align}
The line conductance and susceptance matrices $G$ and $B$ are here assumed to be constant. The specific parameter values for these matrices can be found in~\cite{Zimmerman2011}. In the following, we are additionally interested in the current $I_{i,j}(v,\theta)$ at the transmission line $(i,j) \in \mathcal E$, which can be worked out explicitly in dependence on the voltages $v_j$ and phase shifts~$\theta_j$,
\begin{align}
I_{ij}(v,\theta) = \sqrt{\frac{P_{ij}(v,\theta)^2 + Q_{ij}(v,\theta)^2}{v_i^2}} \; , \notag
\end{align}
where $P_{ij}(v,\theta)$ and $Q_{ij}(v,\theta)$ denote the active and reactive power in the transmission line,
\begin{align}
P_{ij}(v,\theta) \;=\; & v_i^2 G_{ij}  - v_i v_j \left[ G_{ij}  \cos( \theta_i- \theta_j) + B_{ij}\sin(\theta_i- \theta_j) \right] \notag \\
Q_{ij}(v,\theta) \;=\; & v_i v_j \left[ G_{ij}  \sin( \theta_i- \theta_j) + B_{ij}\cos(\theta_i- \theta_j) \right] - v_i^2 B_{ij} \notag
\end{align}
for all $(i,j) \in \mathcal E^+$. The above model equation will next be used to formulate an AC power system state estimation problem.

\subsection{Power system state estimation}
Let $\hat v_i, \hat \theta_i, \hat p_i$, and $\hat q_i$, denote measurements of the voltages, phase shifts and active- and reactive powers at the nodes and let $\hat I_{ij}$ denote measurements for the currents in the transmission lines $(i,j) \in \mathcal E^+$. Next, we consider the least-squares state estimation problem
\[
\min_{v,\theta} \;\;
\sum_{i \in \mathcal N} \left[ \left\|
\left(
\begin{array}{c}
v_i-\hat v_i \\
\theta_i-\hat \theta_i \\
p_i(v,\theta)-\hat p_i \\
q_i(v,\theta)-\hat q_i \\
\end{array}
\right)
\right\|_{\Sigma_i}^2
\hspace{-0.3cm} +\sum_{j\in\mathcal{N}_i}\left\| I_{ij}(v,\theta)-\hat{I}_{ij}\right\|_{\Theta_i}^2
\right].
\]
This problem is not yet in standard form, but it can easily be written in the form~\eqref{eq::opt} by introducing auxiliary variables, where each node keeps copies of the voltage and phase shifts of its neighbors.  In order to set up a realistic case study, the measurements are obtained by running a realistic scenario simulation in~\texttt{MATPOWER} adding randomly generated process noise. Notice that the details of this problem formulation, including the details about how to introduce auxiliary variables and consensus constraints, as well as values for the weighting matrices $\Sigma_i$ and $\Theta_i$ can all be found in~\cite{Du2019}.

\subsection{Numerical results}
\begin{figure}[h]
	\centering	
	\includegraphics[scale=0.2]{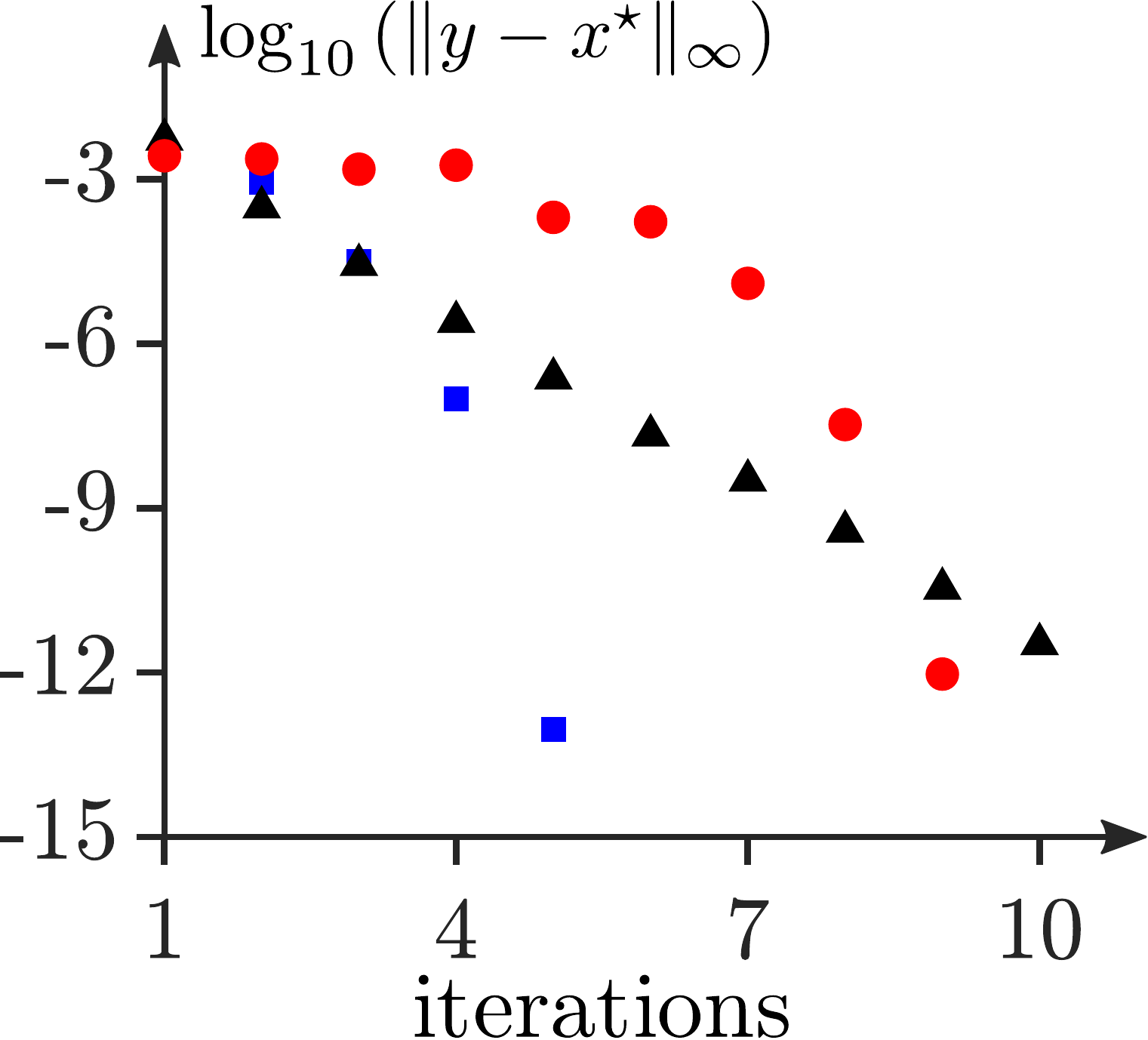}
	\caption{Distance of current iterates to the optimal solution: \texttt{Ipopt} (red circles), Algorithm~\ref{alg::SMM} with exact Hessians (blue diamonds), and Algorithm~\ref{alg::SMM} with Gauss-Newton Hessians (black crosses).}
	\label{fig::case2}
\end{figure}
The diamonds in Figure~\ref{fig::case2} show the numerical results that were obtained by running Algorithm~\ref{alg::SMM} for the above IEEE $33$-bus benchmark case study using models $W_{i,j}$ that are locally accurate up to order $2$, as elaborated in Example~\ref{ex::cubicModel}. The black crosses show results for the same algorithm, but with Gauss-Newton Hessian approximations instead of exact Hessians. These results must be compared to the red circles, which show the iterates of the centralized solver \texttt{Ipopt}~\cite{Wachter2006}. Here, it should be mentioned that, in this figure, one "iteration" of Algorithm~\ref{alg::SMM} refers to a full backward-forward sweep. 
As predicted by Theorem~\ref{thm::local}, we can observe either linear or quadratic local convergence rates of Algorithm~\ref{alg::SMM} depending on the accuracy of the communicated model functions.

\section{Conclusions}
This paper has presented a novel multi-sweep algorithm for asynchronous decentralized optimization over networks with tree graph structure. A first prototype for this method has been presented in the form of Algorithm~\ref{alg::MDP}, which has then been developed further arriving at a fully decentralized multi-sweep communication protocol for non-convex optimization, as presented in Algorithm~\ref{alg::SMM}. We have established conditions under which the proposed method has locally quadratic convergence rate, which have been summarized in Theorem~\ref{thm::local}. This theoretical result has been confirmed numerically by testing the method on a radial AC power network.

\begin{acknowledgements}
	YJ, HY, and BH acknowledge support by ShanghaiTech
	University, Grant-Nr. F-0203-14-012. DK and MD acknowledge support by BMWi via eco4wind (0324125B) and DyConPV (0324166B), and by DFG via Research Unit FOR 2401. 
\end{acknowledgements}

\appendix  %This command ends the counting of sections.
\section{Proof of Theorem~\ref{thm::AlmostNewton}}
\label{app::proof}
Let us introduce the shorthands
\[
z^{k+1} = \left(
\begin{array}{c}
z_1^{k+1} \\[0.16cm]
z_2^{k+1}
\end{array}
\right)  = \left(
\begin{array}{c}
x^{k+1} \\[0.16cm]
\kappa^{k+1}
\end{array}
\right) \qquad \text{and} \qquad 
z^\star = \left(
\begin{array}{c}
x^\star \\[0.16cm]
\kappa^\star
\end{array}
\right)
\]
to denote, respectively, the primal dual minimizer of~\eqref{eq::AlmostNewton} at the $k$-th iteration of the algorithm and the primal-dual minimizer of~\eqref{eq::FC}. Due to the regularity of $x^\star$ the LICQ condition must be satisfied in a neighborhood of $x^\star$, which implies that the first order necessary KKT conditions
\begin{eqnarray}
\label{appeq::aux1}
R(x^k,z^{k+1}) = 0 \qquad \text{and} \qquad R(x^\star,z^\star) = \widetilde R(z^\star) = 0
\end{eqnarray}
with shorthands
\[
R( \xi, \zeta ) = \nabla_z \left[ \Phi(\xi,\zeta_1) + \zeta_2^\tr C(\zeta_1) \right] \qquad \text{and} \qquad \widetilde R(\zeta) = R(\zeta_1,\zeta) = \nabla_z \left[ F(\zeta_1) + \zeta_2^\tr C(\zeta_1) \right]
\]
are satisfied recalling that $\Phi$ is a locally accurate approximation of $F$. Now, because the derivative of $R$ with respect to its second argument, $\nabla_z R(x,\cdot)$, is uniformly Lipschitz continuous function in a neighborhood of $z^\star$, the first equation in~\eqref{appeq::aux1} yields
\begin{eqnarray}
0 &=& R(x^k,z^{k+1}) = R ( x^k, z^k ) + \int_{0}^1 \nabla_z R(x^k, z^k + s (z^{k+1}-z^k) ) (z^{k+1}-z^k) \, \mathrm{d}s \\[0.16cm]
&=& \widetilde R ( z^k ) + M(z_k)  (z^{k+1}-z^k) + \mathbf{O}\left(  \Vert z^{k+1}-z^k \Vert^2 \right) \; ,
\end{eqnarray}
where we have set $M(z^k) = \nabla_z R(x^k, z^k) = \nabla_z \widetilde R(z^k)$ and used that $\widetilde R(z^k) = R(x^k,z^k)$. Notice that the KKT matrix $M(z_k)$ is invertible for all $z^k$ in an open neighborhood of $z^\star$ as we assume that the LICQ and SOSC condition are satisfied at $z^\star$. Consequently, because we have $\widetilde R(z^k) = \mathbf{O}( \Vert z^k - z^\star \Vert )$, the above equation implies that
\begin{eqnarray}
\label{appeq::Newton}
z^{k+1} = z^k - M(z^k)^{-1} \widetilde R(z^k) + \mathbf{O}( \Vert z^k - z^\star \Vert^2 ) \; .
\end{eqnarray}
From here on, the proof is very similar to the standard proof of quadratic convergence of Newton's method (see, e.g.~\cite[Thm.~3.5]{Nocedal2006}); that is we use~\eqref{appeq::Newton} to establish the inequality
\begin{eqnarray}
\Vert z^{k+1} - z^\star \Vert &=& \left\| z^k - z^\star - M(z^k)^{-1} \widetilde R(z^k) \right\| + \mathbf{O}( \Vert z^k - z^\star \Vert^2 ) \notag \\[0.16cm]
&=& \left\| z^k - z^\star - M(z^k)^{-1} \left( \widetilde R(z^k) - \widetilde R(z^\star)  \right) \right\| + \mathbf{O}( \Vert z^k - z^\star \Vert^2 )\notag  \\[0.16cm]
&=& \left\| \left(I - M(z^k)^{-1} \int_0^{1} \nabla_z \widetilde R(z^k+s(z^k-z^\star)) \, \mathrm{d}s \right) (z^k-z^\star)  \right\| \notag \\[0.1cm]
& & + \mathbf{O}( \Vert z^k - z^\star \Vert^2 ) \notag \\[0.16cm]
&=& \underbrace{\left\|  I - M(z^k)^{-1} \nabla_z \widetilde R(z^k)  \right\|}_{=0} \Vert z^k - z^\star \Vert + \mathbf{O}( \Vert z^k - z^\star \Vert^2 ) \; .
\end{eqnarray}
Because the LICQ condition holds the iterates of the multiplier sequence $\kappa^k$ is uniquely determined by the sequence $x^k$  (since $x^{k+1}$ depends only on $x^k$, but not on $\kappa^k$), the above equation also implies that
\[
\Vert x^{k+1} - x^\star \Vert = \mathbf{O}( \Vert x^k - x^\star \Vert^2 ) \; .
\]
The latter equation corresponds to the statement of the theorem establishing local quadratic convergence.

%References
% BibTeX users  please use  
\bibliographystyle{spmpsci_unsrt}
\bibliography{MultiSweep}

\end{document}